\documentclass[12pt,leqno,twoside]{amsart}
\usepackage{amssymb,amsfonts,amsmath,amsthm}
\usepackage{latexsym}
\usepackage{verbatim}
\usepackage{epsfig}
\usepackage[T1]{fontenc}
\usepackage[latin1]{inputenc}

\newtheorem{theorem}{Theorem}[section]

\newtheorem{lemma}[theorem]{Lemma}
\newtheorem{proposition}[theorem]{Proposition}
\theoremstyle{remark}
\newtheorem{remark}[theorem]{\sc Remark}
\theoremstyle{remark}

\theoremstyle{definition}
\newtheorem{definition}[theorem]{Definition}
\theoremstyle{remark}
\newtheorem{example}[theorem]{\sc Example}

\theoremstyle{remark}

\theoremstyle{remark}

\numberwithin{equation}{section}

\renewcommand{\Box}{_\square}    

\newcommand{\cal}{\mathcal}

\newcommand{\reg}{{\rm{reg}}}

\newcommand{\Sing}{\mathop{{\rm{Sing}}}\nolimits}

\newcommand{\mult}{{\rm{mult}}}

\newcommand{\cl}{{\rm{closure}}}

\newcommand{\ity}{{\infty}}

\renewcommand{\d}{{\rm{d}}}

\newcommand{\e}{\varepsilon}

\newcommand{\fin}{\hspace*{\fill}$\Box$\vspace*{2mm}}
\newcommand{\jac}{{\rm{Jac} \ }}

\newcommand{\cB}{{\cal B}}

\newcommand{\cH}{{\cal H}}

\newcommand{\cF}{{\cal F}}

\newcommand{\cS}{{\cal S}}

\newcommand{\cW}{{\cal W}}

\newcommand{\bR}{{\mathbb R}}
\newcommand{\bC}{{\mathbb C}}

\newcommand{\bP}{{\mathbb P}}

\newcommand{\bn}{{\bf n}}

\begin{document}

\title[Curvature and Gauss-Bonnet defect]
 {Curvature and Gauss-Bonnet defect of global affine hypersurfaces}

\author{\sc Dirk Siersma}

\address{Mathematisch Instituut, Universiteit Utrecht, PO
Box 80010, \ 3508 TA Utrecht
 The Netherlands.}

\email{siersma@math.uu.nl}

\author{\sc Mihai Tib\u ar}

\address{Math\' ematiques, UMR 8524 CNRS,
Universit\'e des Sciences et Technologies de  Lille (USTL), \  59655 Villeneuve d'Ascq, France.}

\email{tibar@math.univ-lille1.fr}







\begin{abstract}
 The total curvature of complex
 hypersurfaces in $\bC^{n+1}$ and its variation in families appear to depend not only
 on singularities but also on the 
behaviour in the neighbourhood of infinity. We find the asymptotic
loss of total curvature towards infinity and we  express
the total curvature and the Gauss-Bonnet defect in terms of
 singularities and tangencies at infinity. 
\end{abstract}

\maketitle

\setcounter{section}{0}


\section{Introduction} \label{intro}

Let $Y\subset \bC^{n+1}$ be a global algebraic hypersurface, the zero locus of a
 polynomial in $n+1$ complex variables.
By curvature, denoted by $K$, we mean the Lipschitz-Killing curvature
of a real codimension two analytic space $Y$, 
with respect to the metric induced by the flat Euclidean metric
of $\bC^{n+1}$. Let $\d v$ denote the associated volume form. The
integral of the curvature $\int_{Y}  K \d v$ will be called ``total
curvature'' of $Y$. 

We study here the influence of the position of $Y$ at
 infinity upon the total curvature of $Y$. Computing the 
 total curvature of the  projectivised
 $\bar Y$ wouldn't help, since the metrics on  $\bP^{n+1}$ and $\bC^{n+1}$
are different.
 We shall therefore exploit two ways of computing the total curvature of $Y$: 
(1). by comparing it with the Euler characteristic
 $\chi(Y)$, and (2). by comparing it to the total curvature of a 
general hypersurface, after 
embedding $Y$ into a family.
 
The first approach goes back to extrinsic proofs of the Gauss-Bonnet theorem. 
 The failure of this celebrated theorem in case of
{\em open surfaces} is a theme which has been under constant attention 
ever since Cohn-Vossen's
 pioneering  work \cite{Co} in 1935.  Since our space $Y$ is not 
compact, and possibly singular,  we consider the {\em Gauss-Bonnet defect}:
\[
GB(Y) := \omega_n^{-1}\int_Y K \d v - \chi(Y),
\]
where $\omega_n$ is a universal constant, see \S \ref{ss:real}. 

The second approach is based on the work of Langevin \cite{La-cmh,
  La-thesis} and Griffiths
\cite{Gr} in the late 70's on
 the influence of an isolated singularity upon the total curvature of the
  local Milnor fibre in case of analytic hypersurface germs. 
 Langevin
 found the ``loss of total curvature'' of the Milnor fibre
 at an isolated hypersurface singularity\footnote{More about this topic 
can be found in  Griffiths' paper \cite{Gr}.} and expressed it in terms of 
certain Milnor-Teissier numbers
 $\mu^*$, see \S \ref{main}.

In our global case, we start from the interplay between the total
curvature and the {\em affine class} of $Y$, defined as the number of tangent
hyperplanes to $Y_\reg$ in a general global affine pencil of hyperplanes in
$\bC^{n+1}$. 
 For such pencils one
defines global polar loci which are affine curves. 
We point out here that any global pencil of affine
 hypersurfaces,
even if general, has as a ``limit'' the hyperplane at infinity
$H^{\infty}$, which may be not in general position with respect
to the projectivised hypersurface $\bar Y$.

 We show in Theorem \ref{t:0} that there is a second possible
way of losing curvature when specialising in some family of affine hypersurfaces: 
 towards infinity. The formula for the total curvature can be interpreted as a 
 Pl\" ucker-type formula for the class of {\em affine hypersurfaces}, 
see \S \ref{plucker}. 
In order to get more grip on the meaning of the quantity  of curvature 
absorbed at infinity, we release generality:
 we consider $Y$ with isolated singularities and such that
 $\bar Y \cap H^\ity$ has singularities of dimension $\le 1$. This includes
 the most studied cases in the literature, see \S \ref{ss:nat}. 
We then express the total curvature, as well as the Gauss-Bonnet defect,
 in terms of invariants associated to singularities of $\bar Y$ 
and to the non-generic section $\bar Y\cap H^{\infty}$ of $\bar Y$ (Proposition
\ref{pr:B-type}).
We discuss in \S \ref{examples} 
several examples of deformations of  affine
hypersurfaces $Y$ with isolated and also non-isolated
singularities.


\section{Background on the total curvature}\label{curvature}
\subsection{Real submanifolds}\label{ss:real}
For a real orientable hypersurface of $\bR^{N}$ one has a well-defined Gauss map. One
defines the Gauss-Kronecker curvature $K(x)$ as the Jacobian of the Gauss map
at $x$. For a
submanifold $V$ in $\bR^N$
 Fenchel \cite{Fe} computes the curvature as follows. For a given point $x$ on $V$ one
considers a unit normal vector $\bn$, projects $V$ orthogonally to the affine subspace $W$
generated by the affine tangent space to $V$ and this normal vector. The projection
of $V$ to $W$ is a hypersurface, which has a well defined Gauss-Kronecker curvature $K(x,\bn)$.
The Lipschitz-Killing
 curvature $K(x)$ of $V$ in $x$ is defined (see e.g.\cite[p. 246-247]{ChL})
as the integral of these curvatures over all normal directions,
up to a universal constant $u$:
$K(x) = u \int_{N_xV} K(x,\bn)\d \bn$.

The classical {\em Gauss-Bonnet theorem} says that if $V$ is compact and of
even dimension $2n$ then the
total curvature is equal, modulo an universal constant, to the Euler characteristic:
\[ \omega_n^{-1} \int_V K \d v =  \chi(V),\]
 where $\d v$ denotes the restriction of
the canonical volume form and where $\omega_{n} =
 \frac{(2\pi)^{n}}{1\cdot3\cdots (2n-1)}$ is half the
volume of the sphere $S^{2n}$.
\subsection{Complex hypersurfaces}\label{ss:complex}

Langevin  \cite{La-cmh, La-thesis} studied the integral of curvature of complex
 hypersurfaces $Y\subset
\bC^{n+1}$, using Milnor's approach \cite{Mi} to the
 computation of the total curvature from the number of critical points of
 orthogonal projections on generic lines. We recall here some results and fix our notations.

The curvature $K(x)$ of a smooth complex hypersurface
is the Lipschitz-Killing curvature of $Y_\reg$ as a codimension 2 submanifold of
$\bR^{2n+2}$, where $Y_\reg$ denotes the regular part of $Y$.
A computation due to Milnor allows one to express the Lipschitz-Killing
 curvature of $Y$
in terms of the complex Gauss map $\nu_{\bC} :  Y_\reg \rightarrow \bP^{n}_{\bC}$
which sends a point $x \in Y_\reg$ to the complex tangent space of $Y_\reg$ at $x$, cf
 \cite[pag. 11]{La-cmh}:

\begin{equation}\label{eq:jac}
   (-1)^n  K(x) = | K(x)| = \frac{2 \cdot 4 \cdots 2n}{1\cdot3\cdots (2n-1)} |
 \jac \nu_{\bC} |^2.
\end{equation}
 In the complex case the curvature $K$ is well-known to have the constant
sign $(-1)^n$.
Using (\ref{eq:jac}) one can prove an {\em exchange formula}, as
follows.\footnote{The {\em exchange principle} was originally used
 in the framework
of total absolute curvature of knots and embedded real manifolds,
by Milnor \cite{Mi}, Chern-Lashof \cite{ChL}, Kuiper \cite{Ku}.}
  Let $H$ be a hyperplane in $\bP^{n}$,
defined by a linear form $l_H : \bC^{n+1} \to \bC$.
For almost all $H\in \check \bP^n$
the restriction  of $l_H$ to $Y_\reg$ has only complex Morse critical points.
Let $\alpha_Y(l_H)$ be the number of those critical points (which is finite,
since $Y$ is algebraic).
On the complement of the zero-set of its Jacobian, the complex Gauss map  is a local
diffeomorphism with locally
constant degree $\alpha_Y(l_H)$.  It is shown in \S \ref{s:polar}
 that this
number does not depend on $H$ running in
some Zariski open set of  $\check \bP^n$, the dual projective space of all hyperplanes
of $\bP^{n}$. So we may denote
it by $\alpha_Y$.  From the
above discussion and from Langevin's \cite[Theorem A.III.3]{La-thesis}, one may
draw the following result:
\begin{lemma}\rm (Langevin) \it \label{l:exch}
Let $Y\subset \bC^{n+1}$ be any affine hypersurface. Then:
\[
\int_{Y} | K |\d v  =   \frac{2 \cdot 4 \cdots 2n}{1\cdot3\cdots (2n-1)}
\int_{\check \bP^n} \alpha_Y(l_H) \ \d H = \omega_n \alpha_Y.
\]
\fin
\end{lemma}
Here the integral $\int_{Y} | K |\d v$ is by definition
 the integral over $Y_\reg$. This makes sense since $Y$ differs from $Y_\reg$
by a set of measure zero.
The above formula shows in particular that, up to the constant $\omega_n$,
 $\int_{Y} |K| \d v$ is a non-negative integer.
The real version of the exchange principle can be used to give an extrinsic
 proof of the Gauss-Bonnet theorem for compact even dimensional manifolds.

 In order to measure the failure of the Gauss-Bonnet theorem in case of
 singular or
 non-compact spaces, we use the {\em Gauss-Bonnet defect} of $Y$ defined in
 the Introduction: $GB(Y) := \omega_n^{-1}\int_Y K \d v - \chi(Y)$.
By the above, the Gauss-Bonnet defect of a complex affine hypersurface $Y$ is
an integer. It may be interpreted
as the correction term due to the ``boundary at infinity'' of $Y$, at least in
case $Y$ has isolated singularities, as follows.

Let $B_R\subset \bC^{n+1}$ be a ball centered at the origin and denote 
$Y_R := Y \cap B_R$ and
$\partial Y_R := Y\cap \partial \bar B_R$. Since $Y$ has isolated
singularities and is affine, the intersection
$Y\cap \partial \bar B_R$ is transversal and $Y_R$ is diffeomorphic to $Y$,
 for large enough radius $R$.
By applying the
Gauss-Bonnet formula for the manifold with boundary $Y_R$, 
see Griffith \cite[p. 479]{Gr}, we get:
\[  \omega_n^{-1}\int_Y K \d v - c \int_{\partial Y_R} k \d s = \chi(Y_R),\]
where $k$ is the generalised `geodesic curvature' of $\partial Y_R$ and $c$ is
a universal
constant (which will be not made precise here). It then follows:
\[ GB(Y)= \lim_{R\to \infty} c \int_{\partial Y_R} k \d s.\]
This interpretation suggests that $GB(Y)$ should be related to singularities
which occur at infinity. We shall find such a relation in (\ref{eq:B_1-form}). 
\subsection{Pl\" ucker's class formula}\label{ss:class}
Let $V \subset \bP^{n+1}$ be a projective hypersurface of degree $d$.
The space of tangent hyperplanes 
to $V_\reg$ is a subset in the dual 
$\check \bP^{n+1}$ and its closure $\check{V}$ is called the dual of $V$.
The degree of $\check{V}$, denoted by $d^*(V)$, is the number of intersection 
points of $\check V$ with a generic
projective line in $\check \bP^{n+1}$. This is the same as the {\em class} of $V$,
 the number of tangent hyperplanes to $V_\reg$ in a generic pencil on $\bP^{n+1}$.
Pl\" ucker's class formula describe $d^*$ in terms of $d$ and of certain
invariants of the singularities of $V$.
The one proven by Pl\" ucker himself in 1834
considers curves with nodes and cusps. Teissier generalized it in 1975 to the 
case of projective hypersurfaces
with isolated singularities, and  Laumon \cite{Lau} found the following
 equivalent formula, in terms of Milnor-Teissier numbers of isolated
 singularities
(see \S \ref{sss:polarisol}):
\begin{equation}\label{eq:class}
  d^*(V)  = d(d-1)^{n} -
\sum [\mu^{\langle n\rangle} + \mu^{\langle n-1\rangle}].
\end{equation}

Later Langevin \cite{La-rendi, La-thesis} showed the
connection with the complex Gauss map and provided
 the integral-geometric interpretation of (\ref{eq:class}). 
Further generalisations, for arbitrary
{\em projective} varieties with isolated singularities, and then without conditions on
singularities, were found notably by
Kleiman, Pohl and respectively Thorup, see e.g. \cite{Th}.

Turning now back to the affine case: the positive integer $\alpha_Y(l_H)$ defined
at \S \ref{ss:complex} can be interpreted 
as the degree of the dual variety $\check Y$. We shall derive in \S \ref{plucker}
 a formula for the affine class of $Y$.

\section{Polar invariants, singularities and Euler characteristic}\label{polar}

The use of polar methods is naturally suggested by the exchange principle.
 On the other hand, 
the affine polar invariants determine, via the Lefschetz slicing theory,
 a CW-complex structure of the space, and therefore
its Euler characteristic.
\subsection{Polar curves in affine families, after \cite{Ti-compo}}
Let $\{X_s\}_{s\in \delta}$ be a family of affine hypersurfaces
$X_s \subset \bC^{n+1}$, where $\delta$ is a small disk at the origin of
$\bC$.
We assume that the family is polynomial, i.e. there is a polynomial
$F :\bC \times \bC^{n+1} \to \bC$ such that $X_s = \{ x\in \bC^{n+1}\mid
F_s(x)= F(s,x)=0\}$. Let us
denote by $X = \cup_{s\in \delta} X_s$ the total space of the family, which is
itself a hypersurface in $\delta \times  \bC^{n+1}$.
Let $\sigma : X \to \delta \subset \bC$ denote the projection of $X$
to the first factor of $\bC \times \bC^{n+1}$.

Let our affine hypersurface $X\subset \bC\times\bC^{n+1}$ be stratified by its canonical (minimal)
Whitney stratification $\cS$, cf. \cite{Te2}. This is a finite stratification, having $X\setminus
\Sing X$ as a stratum. For instance, if $X_s$ has no singularities and $X_0$ has at most isolated
ones, then $\cS$ has as lower dimensional strata only these singular points.
  We shall use the same notation  $l_H$ for
the application $\bC \times  \bC^{n+1} \to \bC$, $(s,x) \mapsto l_H(x)$, as
well as for its restriction to $X$.
The
{\em polar locus} of the map $(l_H,\sigma) : X \to \bC^2$ with respect to $\cS$ is the
following analytic set:
\[ \Gamma_\cS (l_H,\sigma) := \cl \{ \Sing_\cS (l_H,\sigma) \setminus
(\Sing_\cS l_H \cup \Sing_\cS \sigma)\} ,\] where $\Sing_\cS \sigma := \bigcup_{\cS_i \in \cS}
\Sing \sigma_{| \cS_i}$ is the singular locus of $\sigma$  with respect to $\cS$. The singular
loci $\Sing_\cS l_H$ and $\Sing_\cS (l_H,\sigma)$ are similarly defined.


\begin{lemma}\label{l:bertini}{\rm \cite{Ti-compo}} \
There is a Zariski-open set $\Omega_{\sigma} \subset \check\bP^{n}$
such
that, for any $H\in \Omega_\sigma$, the polar locus $\Gamma_\cS (l_H,\sigma)$
is a
curve or
it is empty.
\hspace*{\fill}$\Box$
\end{lemma}
Let $\Omega_\sigma$ be the Zariski-open set from Lemma \ref{l:bertini}. We
 denote by
$\Omega_{\sigma,0}$ the Zariski-open set of hyperplanes $H\in \Omega_\sigma$
 which are transversal
to the canonical Whitney stratification of the projective hypersurface 
$\overline{X_0} \subset
\bP^{n+1}$. This supplementary condition insures that
$\dim(\Gamma_\cS(l_H,\sigma)
 \cap X_0) \le
0$, $\forall H \in \Omega_{\sigma,0}$.
\subsection{The $\alpha^*$ sequence and the Euler characteristic}\label{s:polar}
Let $\{X_s\}_{s\in \delta}$  be any family as above.
We have defined in \cite[\S 3]{Ti-equi} generic polar intersection
multiplicities for such a family. We shall paraphrase that definition by
considering only the regular part  $(X_s)_\reg$ of the hypersurfaces, as
follows:
\begin{definition}\label{d:polarint}
Let $H \in \Omega_{\sigma,0}$.
The following global generic polar intersection multiplicity:
\begin{equation}\label{eq:alpha}
\alpha_{X_s}^{(n)} = \mult (\Gamma_\cS(l_H, \sigma), (X_s)_\reg).
\end{equation} 
is well defined for any $s\in \delta$ and does not depend on the choice
 of $H\in \Omega_{\sigma,0}$.
\end{definition}
The  geometric interpretation of $\alpha_{X_s}^{(n)}$  is the number of Morse 
points of a generic linear function
on $(X_s)_\reg$. 
We shall next define the lower global polar intersection multiplicities 
$\alpha_{X_s}^{(i)}$ by following \cite[\S 3]{Ti-equi}. The idea is to
consider successively general hyperplane slices of our family and apply
Definition \ref{d:polarint}. This idea comes from Teissier's 
construction of polar multiplicities \cite{Te0, Te1, Te2}. 

 One takes a general hyperplane
 $\cH\in \Omega_{\sigma,0}$ and denotes by
$\alpha^{(n-1)}_{X_s}$ the global generic polar intersection multiplicity
 at $s\in \delta$ of the
family of affine hypersurfaces $X' =  X\cap \cH$. 
One pursues in this way and defines step-by-step
$\alpha_{X_s}^{(n-i)}$,
 for $1\le i\le n-1$. We set $\alpha_{X_s}^{(0)} := \deg X_s$.

By a standard connectivity argument, the polar intersection
multiplicities $\alpha_{X_s}^{(i)}$ do not depend on the choices of generic
hyperplanes.
They are also invariant up to linear changes of
 coordinates
 but not invariant up to nonlinear changes of coordinates
 (e.g. $\deg X_s$ is not invariant).
 The numbers $\alpha_{X_s}^{(i)}$ are constant on
 $\delta \setminus \{ 0\}$, provided that
 $\delta$ is small enough. 

The geometric
interpretation of the sequence of global generic polar multiplicities  
$\alpha^{(i)}$ also follows, as we have shown above.
For instance, if we apply this construction to a single
non-singular hypersurface $Y\subset \bC^{n+1}$, by the Lefschetz slicing
principle we get
that $Y$ has the structure of a  CW-complex of
dimension $\le n$, with 
$\alpha^{(i)}_Y$ cells in dimension $i$. 
Consequently, one may expressed its 
Euler characteristic as follows:
\begin{equation}\label{eq:eulerchar}
\chi(Y) = \sum_{i=0}^n (-1)^{i} \alpha^{(i)}_Y.
\end{equation}
In case of a singular $Y$, the formula needs correction; we explain here the case
 of isolated singularities, which we shall use in \S \ref{curv} (and send to
 \S \ref{ss:exnonisol} for non-isolated 
singularities and several examples). 
 By the stratified Morse theory \cite{GM} and Lefschetz slicing principle, 
the space $Y$ is obtained
from the generic slice $Y\cap \cH$ by attaching cones over the complex links 
of each Morse stratified
singularity of the generic pencil on $Y$. The singularities of the pencil on $Y_\reg$
contribute by $\alpha^{(n)}_Y$. In case $Y$ has only isolated singularities, 
the
contribution at each such point-stratum is precisely  the
Milnor number of the generic local hyperplane section  (since our pencil is locally
 generic at those points), which, by a standard argument, is equal to the
 sectional Milnor-Teissier number $\mu^{\langle
 n-1\rangle}$ (see after (\ref{eq:langevin}) for the notation). The slice
 $Y\cap \cH$ and the lower dimensional ones are non-singular.
 We therefore get, in case $Y$ has isolated singularities, the following formula:
\begin{equation}\label{eq:beta2}
\chi(Y) = \sum_{i=0}^n (-1)^{i} \alpha^{(i)}_Y + (-1)^n \sum_{q\in \Sing Y} 
\mu_q^{\langle n-1 \rangle}(Y).
\end{equation} 

\section{Vanishing curvature and an affine Pl\" ucker formula}\label{main}

\subsection{The vanishing curvature}\label{ss:loss}
We show here that in case of a
 family of affine hypersurfaces,  part
of the ``loss of total curvature'' may occur at infinity. We shall denote by 
$\complement  B_R$ the complement in $\bC^{n+1}$ of the ball $B_R$ centered at
 the origin and of radius $R$.
We shall use the shorter notation $\alpha_{s}^{(n)}$ for $\alpha_{X_s}^{(n)}$ 
in the rest of
 the paper.
\begin{theorem}\label{t:0}
 Let $Y\subset \bC^{n+1}$ be any hypersurface. Let
$\{ X_s \}_{s\in \delta}$ be a
one-parameter deformation of $X_0 := Y$ such that $X_s$
 is non-singular for all $s\not=0$. Then the following limit exists:
\begin{equation}\label{eq:main}
\lim_{s\to 0} \omega_{n}^{-1}\int_{X_s} | K| \d v
=  \omega_{n}^{-1}\int_{X_0} | K| \d v  + 
\mult(\Gamma_\cS(\sigma, l_H),
X_0) + \alpha_0^{(n)}(\ity),
\end{equation}
where $\alpha_0^{(n)}(\ity)$ is a non-negative integer defined as:
\begin{equation}\label{eq:main2}
\alpha_0^{(n)}(\ity) := \omega_{n}^{-1} \lim_{R\to \ity}\lim_{s\to 0}\int_{X_s
  \cap \complement  B_R} | K| \d v. 
\end{equation}
\end{theorem}
\begin{proof}
We deduce from Lemma \ref{l:exch} the following
 general formula, by using Definition (\ref{d:polarint}):
\begin{equation}\label{eq:intgen}
\omega_n^{-1} \int_{X_s} K \d v = (-1)^n\alpha_s^{(n)}.
\end{equation}
It has been remarked in \S \ref{s:polar} that  $\alpha_s^{(n)}$ is constant for 
$s\in \delta \setminus \{ 0\}$, if the disk $\delta$ is small enough.
Therefore the limit $\lim_{s\to 0} \omega_{n}^{-1}\int_{X_s} | K| \d v$
 is equal to $\alpha_s^{(n)}$.

Let us take  $H\in \Omega_{\sigma,0}$ as in \S \ref{polar}. From
the definition (\ref{eq:alpha}) of $\alpha_s^{(n)}$ we get the following
decomposition into a sum of intersection numbers:
\begin{equation}\label{eq:sum}
\alpha_s^{(n)} = \alpha_0^{(n)} + \alpha_0^{(n)}(crt, H) +
\alpha_0^{(n)}(\infty, H).
\end{equation}
The first term is the intersection multiplicity 
$\mult (\Gamma_\cS(l_H, \sigma), (X_0)_\reg)$ and we know that it does not depend
on the choice of $H$ as above and that it is equal to 
$\omega_{n}^{-1}\int_{X_0} | K| \d v$, which is the first term in our claimed
formula (\ref{eq:main}).
The second term of the sum (\ref{eq:sum}) counts the number of
those intersection points of $\Gamma_\cS(\sigma, l_H)$ with $(X_s)_\reg$ which tend
to points $q\in \Sing X_0$. This multiplicity does
not depend on the choice of generic $H$.
It then follows that the third term from (\ref{eq:sum}), namely
$\alpha_0^{(n)}(\infty, H)$,
is also independent on $H\in\Omega_{\sigma,0}$. It counts 
 the asymptotic loss of
intersection points of the polar curve $\Gamma_\cS(l_H, \sigma)$ with $(X_s)_\reg$,
as $s\to 0$. In other words, we have:
\begin{equation}\label{eq:ident} 
\alpha_0^{(n)}(\infty, H) = \lim_{R\to \ity}\lim_{s\to 0}
 \mult(\Gamma_\cS(\sigma, l_H), (X_s)_\reg\cap \complement  B_R).
\end{equation}
 Let us see that this is exactly the double limit defined 
by (\ref{eq:main2}). By the exchange formula (Lemma
\ref{l:exch}) we have that:
\[ \int_{X_s\cap \complement B_R} | K| \d v = u\int_{\check \bP^n}
 \alpha_{X_s\cap \complement B_R}(l_H) \ \d H , \]
where $u$ is a constant defined in Lemma
\ref{l:exch}.
Since this integral is, by definition, bounded from above by $\omega_n \alpha_s^{(n)}$,
we may apply Lebesgues's theorem of dominated convergence (also used by
Langevin in his local proof \cite{La-cmh}). This allows us to interchange 
each of the limits with the integral, thus we get:
\[ \lim_{R\to \ity}\lim_{s\to 0}\int_{X_s
  \cap \complement  B_R} | K| \d v = u\lim_{R\to \ity}\lim_{s\to 0} \int_{\check \bP^n}
 \alpha_{X_s\cap \complement B_R}(l_H) \d H =  u\int_{\check \bP^n}
 [\lim_{R\to \ity}\lim_{s\to 0}\alpha_{X_s\cap \complement B_R}(l_H)] \d H.
\]
Since $\alpha_{X_s\cap \complement B_R}(l_H) = \mult(\Gamma_\cS(\sigma, l_H),
(X_s)_\reg\cap \complement  B_R)$, by using now (\ref{eq:ident}) we get our
 claimed equality.
\end{proof}
In case of a nonsingular $X_0$, the non-negative integer
 $\alpha_0^{(n)}(\ity)$ is precisely the
``polar defect at infinity'' which has been
 introduced in \cite{Ti-equi} under the notation $\lambda_0^{n}$.
We shall see in \S \ref{curv} 
how $\alpha_0^{(n)}(\ity)$ can be expressed in terms of singularities occurring
at infinity, in certain situations. 
\subsection{A general Pl\" ucker-type formula for the class of 
affine hypersurfaces}\label{plucker}
Let $Y \subset \bC^{n+1}$ be a hypersurface of degree $d$.
The degree $\deg(\check Y)$ of the affine dual $\check Y$
is equal to the number of tangent
hyperplanes to $Y_\reg$ in a  generic affine pencil of hyperplanes in
$\bC^{n+1}$. We shall call it the {\em affine class} of $Y$ 
in analogy to the
projective case (see \S \ref{ss:class}), and we shall denote it by $d^@ (Y)$.
The affine pencils (see \S \ref{polar}) differ from the projective pencils
especially in a neighbourhood of infinity, since after projectivising, the
hyperplane at infinity $H^\ity$ becomes a member of the pencil and our
hypersurface $Y$ may be asymptotically tangent to $H^\ity$. 

We say that
an affine hypersurface of degree $d$ is {\em general} when
its projective closure is non-singular and transverse to
the hyperplane at infinity. Its Euler characteristic is equal 
to $1 + (-1)^n (d-1)^{n+1}$. The polar intersection number $\alpha^{(n)}$
 (Definition (\ref{eq:alpha})) is then maximal; by B\'ezout theorem, it is equal
 to $d(d-1)^n$. 

Next, let us remark that one may always deform $Y :=X_0$
 in a constant degree family such that $X_s$ is
general, for $s\not= 0$. For instance, for $X_0 := \{f=0\}$, 
define $f_s = (1-s)f + s(g_d
-1)$, where $g_d = x_1^d + \cdots + x_{n+1}^d$. Then  $X_s :=
\{f_s=0\}$ has this property, for small enough $s\not= 0$.

Considering some deformation of $Y=X_0$
in a constant degree family of general hypersurfaces, we may derive
 from Theorem \ref{t:0} the following formula for the affine class:

\begin{equation}\label{eq:affine}
 d^@ (X_0) = d(d-1)^n - \mult(\Gamma_\cS(\sigma, l_H),
X_0) - \alpha_0^{(n)}(\ity).
\end{equation}

\smallskip
\noindent
 This can be made more explicit in case of isolated singularities,
 with help of the forthcoming formulas (\ref{eq:mu-teiss})
and (\ref{eq:alpha2}).



\subsection{Case of isolated affine singularities}\label{ss:isol}
\subsubsection{Polar multiplicity}\label{sss:polarisol}
In our global case, if $X_0=Y$ has only isolated singularities, then one may
identify the intersection multiplicity in the formula (\ref{eq:main}) as
follows:
\begin{equation}\label{eq:mu-teiss}
\mult(\Gamma_\cS(\sigma, l_H),
X_0) = \sum_{q\in \Sing X_0}[\mu_q^{\langle n-1\rangle}(X_0) + \mu_q^{\langle
  n\rangle}(X_0)].
\end{equation}
This comes from the equality for the generic local polar 
multiplicity:
\begin{equation}\label{eq:locmult} 
\mult_q(\Gamma_\cS(\sigma, l_H),
X_0) = \mu_q^{\langle n\rangle}(X_0) + \mu_q^{\langle n-1\rangle}(X_0)
\end{equation}
proved by Teissier \cite{Te1, Te2}  
when $X_0$ is the germ of the zero locus of a holomorphic
function $(\bC^{n+1},0) \to \bC$.
 It is actually well-known that the local equality (\ref{eq:locmult}) is 
valid for any
smoothing of $X_0$. In our case the local smoothing is embedded in 
the global smoothing $\sigma : X \to \bC$.

In the local case, for a {\em germ of
a holomorphic function} with isolated singularities
 $g : (\bC^{n+1},0) \to (\bC,0)$, Langevin's
 formula \cite[Th\' eor\` eme 1]{La-cmh} shows that the loss of total
 curvature at an isolated singularity is measured by the sum of the first two
Milnor numbers of the sequence $\mu^*$ defined by Teissier 
\cite{Te1}:
\begin{equation}\label{eq:langevin}
 \lim_{\e \to 0}\lim_{t\to 0}\int_{g^{-1}(t)\cap B_\e} | K | \d v =
\omega_{n}(\mu^{\langle n\rangle}
+\mu^{\langle n-1\rangle}). 
\end{equation}
 Here $\mu^{\langle n\rangle}$ denotes the
usual Milnor number and $\mu^{\langle n-1\rangle}$ is the Milnor number of a
generic hyperplane section\footnote{Indices are shifted by -1 from
  the original Teissier notation.}.
 The sum
$\mu^{\langle n\rangle} + \mu^{\langle n-1\rangle}$
is precisely the local {\em generic polar number} of $g$, 
i.e. the intersection number of the local polar curve $\Gamma (g,l_H)$
with $g^{-1}(0)$. 

\subsubsection{Gauss-Bonnet defect}\label{sss:GBisol}
For affine $Y$ with isolated singularities (still without any condition 
at infinity), the following formula follows
from (\ref{eq:beta2}) and (\ref{eq:intgen}):
\begin{equation}\label{eq:gb1}
GB(Y) =  (-1)^{n-1} \sum_{q\in \Sing Y}\mu_q^{\langle n-1\rangle}(Y)
- \sum_{i=0}^{n-1}  (-1)^i \alpha_Y^{(i)},
\end{equation}
where the sum $\sum_{i=0}^{n-1}  (-1)^i \alpha_Y^{(i)}$ is just 
$\chi(Y\cap \cH)$.

\section{Total curvature and singularities at infinity}\label{curv}
We focus in the remainder on explaining the loss at infinity of 
the total curvature and the Gauss-Bonnet defect in some distinguished classes of
hypersurfaces.
\subsection{Some natural classes of hypersurfaces}\label{ss:nat}
\begin{definition}\label{d:FB} 
\begin{itemize}
\item[(i)] $Y$  is a  {\em $\cF$-type} hypersurface if $\bar Y$
and $\bar Y\cap H^\ity$
have at most isolated singularities.
\item[(ii)] $Y$  is a {\em $\cB_0$-type}  hypersurface if  $\bar Y$
has at most isolated singularities.\footnote{The topology 
of $\cB_0$-type
  polynomials has been studied in several papers, see e.g. Broughton's
  \cite{Br} and \cite{STexch}.}
\item[(iii)] $Y$  is a {\em $\cB_1$-type}  hypersurface if $Y$ has at most isolated
  singularities
and $\bar Y\cap H^\ity$ has at most isolated 1-dimensional
  singularities.
 \end{itemize}
\end{definition}
It is easy to see that $\cF$-type $\subset$  $\cB_0$-type $\subset$ $\cB_1$-type. 

 In order to introduce the main result of this section, we need to consider
 generic hyperplanes, in the following sense. 
Let $\cW$ be some Whitney stratification of $\bar Y$ such that
 $\bar Y\cap H^\ity$ is a union of strata.
 Let  $\cH \subset \bC^{n+1}$ be a hyperplane such that $\bar \cH$ is 
generic with respect to the strata of $\cW$. There exists a Zariski-open subset of
such hyperplanes, see \S \ref{polar} for a similar discussion. 
\begin{proposition}\label{pr:B-type}
\begin{enumerate}
\item If $Y$ is a $\cB_1$-type hypersurface of degree
$d$ then:

\begin{equation}\label{eq:B_1-form} 
GB(Y)=  (-1)^n (d-1)^n -1  + (-1)^{n+1} \sum_{q\in \Sing Y} 
\mu_q^{\langle n-1\rangle}(Y) + 
 \end{equation} 
\[ 
(-1)^{n+1}[ \mu_p(\bar Y\cap \bar\cH) +\mu(\bar Y\cap\bar\cH \cap
H^\ity)].\]
\smallskip
\item If  $Y$ is a $\cB_0$-type hypersurface of degree
$d$ then:

\begin{equation}\label{eq:B-form}
\omega_n^{-1} \int_{Y} K \d v = (-1)^n d(d-1)^n  + (-1)^{n+1} \sum_{q\in \Sing
  Y} [\mu_q^{\langle n\rangle}(Y) +\mu_q^{\langle n-1\rangle}(Y)] + 
 \end{equation}
\[ (-1)^{n+1}\sum_{p\in (\Sing \bar Y) \cap H^\ity} \mu_p(\bar Y) +
(-1)^{n+1}\mu(\bar Y\cap\bar\cH \cap H^\ity) +
\chi^{n,d} - \chi(\bar Y \cap H^{\ity}),\]
 where $\chi^{n,d}$ denotes the Euler characteristic of the generic
 hypersurface
of degree $d$ in $\bP^n$ and where
 $\mu(\bar Y  \cap \bar\cH \cap H^\ity)$ is a notation for
 $\sum_{p\in \bar \cH\cap
  \Sing (\bar Y \cap H^\ity)} \mu_p(\bar Y \cap \bar\cH \cap H^\ity)$.
\end{enumerate}
\end{proposition}
\smallskip
 
\begin{proof}
We have by (\ref{eq:gb1}): $GB(Y) = 
 (-1)^{n-1} \sum_{q\in \Sing Y}\mu_q^{\langle n-1\rangle}(Y)
- \chi(Y\cap \cH)$. Now $Y\cap \cH$ is $\cF$-type and we may 
compute its Euler characteristic by taking a deformation of $Y=X_0$
in a constant degree family such that $X_s$ is general for $s\not= 0$,
as follows:

\[  \chi(X_s\cap \cH) - \chi(X_0\cap \cH) 
   = [\chi(\bar X_s\cap \bar\cH) - \chi(\bar X_0\cap \bar\cH)] + [
  -\chi(\bar X_s\cap \bar\cH \cap H^\ity) + \chi(\bar X_0\cap \bar\cH \cap
  H^\ity)]\]
\[=
 (-1)^{n-1}\sum_{p\in \Sing (\bar X_0 \cap \bar\cH)}
 \mu_p(\bar X_0\cap \bar\cH) + (-1)^{n-1}
 \sum_{p\in \Sing (\bar X_0 \cap \bar\cH\cap H^\ity)}
 \mu_p(\bar X_0 \cap \bar\cH\cap H^\ity).\]

\smallskip
 \noindent
We then get (\ref{eq:B_1-form}) since $\chi(X_s\cap \cH) = 1 + (-1)^{n-1} (d-1)^n$.
\smallskip
Let us prove (b) now. For the $\cB_0$-type hypersurface $Y= X_0$,
the singularities of $\bar Y$ are isolated but those of $\bar Y \cap
H^\ity$ are of dimension at most 1. We therefore have:

\[
\chi(X_0) -\chi(X_s) =  \chi(\bar X_0) -\chi(\bar X_s) -  \chi(\bar X_0 \cap
H^\ity) +\chi(\bar X_s \cap H^\ity)=
 \] \[
(-1)^{n+1} \sum_{q\in \Sing X_0} \mu_q^{\langle n\rangle}(X_0) +(-1)^{n+1}\sum_{p\in (\Sing \bar X_0) \cap H^\ity}
\mu_p(\bar X_0) + \chi^{n,d} - \chi(\bar X_0 \cap H^\ity).
\]

\smallskip
\noindent
We then get our result from the definition of $GB$, by using the 
 equality (\ref{eq:B_1-form}).
\end{proof}
Comparing (\ref{eq:B-form}) to (\ref{eq:main}) and to 
 (\ref{eq:mu-teiss}) we get, for a deformation of $Y=X_0$
in a constant degree family such that $X_s$ is general for $s\not= 0$:

\begin{equation}\label{eq:alpha2} 
\alpha_0^{(n)}(\ity) =
\sum_{p\in (\Sing \bar Y) \cap H^\ity} \mu_p(\bar Y) +
\mu(\bar Y\cap\bar\cH \cap H^\ity) + (-1)^{n+1}
[\chi^{n,d} - \chi(\bar Y \cap H^{\ity})].
\end{equation}
\begin{remark}\label{r:th}
As a particular 
case of (\ref{eq:B-form}), the following
formula holds for an $\cF$-type hypersurface:
\begin{equation}\label{eq:F-form}
\omega_n^{-1} \int_{Y} |K| \d v  =  d(d-1)^n
-\sum_{q\in \Sing Y} [ \mu_q^{\langle n\rangle}(Y)
+\mu_q^{\langle n-1\rangle}(Y)] 
\end{equation}
\[- \sum_{p\in  \Sing (\bar Y  \cap  H^\ity)}[\mu_p(\bar Y) +
 \mu_p(\bar Y \cap H^\ity) ]. \] 
The contribution from the affine singularities 
is contained in the first of the two sums: one
recognizes the Milnor-Teissier numbers of formula (\ref{eq:class}). 
The second sum is due to the ``singularities at infinity'':
the number 
$\mu_p(\bar Y) +
 \mu_p(\bar Y \cap H^\ity)$ is exactly 
the local polar number $\lambda_p =
 \mult_p(\Gamma (\sigma ,x_0), \bar X_0)$ of
  the polar curve of the family $\{ X_s\}_s$ with respect to the local
coordinate at infinity $x_0$, which is {\em not} a locally generic 
coordinate\footnote{In this context, it
 was used in \cite[3.7]{Ti-equi}.}
 (compare to \S \ref{sss:polarisol}).
 Local polar numbers, introduced 
by Teissier in \cite{Te0}, are well defined as soon as
the polar locus is a curve.

We shall give an example of a $\cF$-type family specialising to 
 a $\cB_0$-type hypersurface, such that
 the  Euler characteristic is constant but the total curvature jumps
 (Example \ref{e:2}).
\end{remark}
\subsection{Concentration of the loss of total curvature at infinity}
\label{c:f-type}
In the case of 
$\cF$-type hypersurfaces
there is pointwise
 concentration of the loss of total curvature at infinity, see  (\ref{eq:F-form}).
  This might be no longer the case for $\cB$-type or more general
 classes of hypersurfaces: in formula (\ref{eq:B-form}) we have
Euler characteristics and dependence on the
 slice $\bar \cH$. 
 The loss of total curvature
 at infinity is nevertheless concentrated at the singular locus of the set $\bar Y \cap
 H^{\ity}$.



\subsection{Affine curves and the correction term at infinity}\label{ss:risler}
 Let $C := \{ f=0\}\subset \bC^2$ be
 a non-singular complex affine curve of degree
$d$. We get from (\ref{eq:B_1-form}):
\begin{equation}\label{eq:c} 
GB(C) = -d. 
\end{equation}
The well-known
 inequality due to Cohn-Vossen \cite{Co} tells that
$GB(M) \le 0$ if $M$ is a complete, finitely connected Riemann surface having
absolutely integrable Gauss curvature. 

 Let now $r$ be the number of asymptotic directions of $C$, i.e. the
number of points in the set $\{ f_d =0\}$, where $f_d$
 denotes the degree $d$ homogeneous part of $f$. 
Let us point out that since $\bar C  \cap
 H^\ity$ consists of $r$ points, the sum of Milnor numbers $\sum_{p\in  
\Sing (\bar C  \cap  H^\ity)}\mu_p(\bar C \cap
H^\ity)$ is precisely $d- r$. By applying formula (\ref{eq:F-form}),
 since non-singular plane curves are of $\cF$-type, we get:
\begin{equation}\label{eq:risler} 
\begin{array}{c}
\omega_n^{-1} \int_{C} |K| \d v  =  d(d-1)
- \sum_{p\in  \Sing (\bar C  \cap  H^\ity)}\mu_p(\bar C) - d +  r = \\

  \ \ \ \ \  \ \ = d^2 -2d +   r -  \sum_{p\in  \Sing (\bar C 
 \cap  H^\ity)}\mu_p(\bar C).
\end{array}
\end{equation}
Comparing this to the formula found by Risler 
 \cite[Proposition 4.2, (15)]{Ri-bullms} for a non-singular complex affine
 curve $C$, one notices that the latter  
does not contain the sum $\sum_{p\in  \Sing
  (\bar C  \cap  H^\ity)}\mu_p(\bar C)$. Therefore Risler's formula
 would not be  valid when
the compactification $\bar C$ is singular. However, Risler uses his
 formula in {\em loc.cit.} only in the case $r=d$, which 
implies that the affine curve $C$ is 
 general at infinity. In this special case indeed formula
 (\ref{eq:risler}) reduces as such.
\subsection{Semi-continuity and extrema of curvature integrals}
For any family $\{X_s\}_{s\in \delta}$ of affine hypersurfaces we have:
\begin{equation}\label{eq:scont1}
 \omega_n^{-1} \int_{X_0} |K| \d v =  \alpha_0^{(n)} \le \alpha_s^{(n)}
 = (\omega_n)^{-1} \int_{X_s} |K| \d v.
\end{equation}
The total curvature is therefore bounded as follows: $0\le 
\omega_n^{-1} \int_{X_0} |K| \d v \le d(d-1)^n$. 

For a general hypersurface $X_0$, the  equality 
$\omega_n^{-1} \int_{X_0} |K| \d v = d(d-1)^n$
 holds. We claim that the reciprocal is true.
Indeed, if $X_0$ is not general then there exists a deformation 
$\{ X_s\}_s$ such
that: $X_s$ is of $\cF$-type for $s\not= 0$, $\bar X_0$ is non-singular
and $\bar X_0 \cap H^\ity$ is  non-singular except at one point, say $p$, where
the singularity is of type $A_1$, i.e. $\mu_p(\bar X_0 \cap
H^\ity)=1$. According to (\ref{eq:F-form}) we then have: 
  $\omega_n^{-1} \int_{X_s} |K| \d v = d(d-1)^n -1$, 
which, together with the semi-continuity relation (\ref{eq:scont1}), 
gives a contradiction.

  What happens now when the minimum occurs, i.e. the total curvature of 
$X_0$ is zero?
 For the case of non-singular $X_0$,
  the answer is the following:
$(\omega_n)^{-1} \int_{X_0} |K| \d v = 0$ implies that the map $l_H :
 X_0 \to \bC$, for $H\in \Omega_\sigma$, is a trivial fibration; 
in particular $b_n(X_0) =0$. 
This is a consequence of the fact that $\alpha_0^{(n)}=0$ implies that there are no
$n$-cells
in the CW model of $X_0$, see \S \ref{s:polar}.

\section{examples}\label{examples}
\begin{example}\label{e:1}
We show first how to compute the total absolute curvature directly from
equation (\ref{eq:intgen}).
Let
$f: \bC^3 \to \bC$, $f(x,y,z) = x+ x^2yz$. We consider the family $X_s = \{
f=s\}$, see \cite[Example 3.8]{Ti-equi}.
The generic polar intersection multiplicities and the
 defects at
 infinity in the neighbourhood of the value $0$ are given in \cite{Ti-equi};
 from those results we may extract the folowing data: $\alpha^{(2)}_s = 5$, 
$\alpha^{(1)}_s = 8$, $\alpha^{(0)}_s = 4$ for $s \ne 0$,
and $\alpha^{(2)}_0 = 3$, $\alpha^{(1)}_0 = 6$, $\alpha^{(0)}_0 = 4$.
 We get: 
$\omega_2^{-1} \int_{X_s} | K |\d v  =  5$ if $s \ne 0$ and
$\omega_2^{-1}\int_{X_0} | K |\d v  =  3$.

The variation of total curvature is $2$ and is equal 
to the vanishing curvature at infinity
$\alpha_0^2(\ity)$, as defined
in  Theorem \ref{t:0}. Therefore the curvature of $X_s$ is not constant in the
family, even if $X_s$ is nonsingular and
$\chi(X_s) =1$ for all $s\in \bC$ (see {\em loc.cit.}).  
It is also clear that the family is not topologically trivial,
since the number of connected components of the fibers change at $s=0$.

\end{example}


\begin{example}\label{e:2}
Consider the double parametre family $X_{s,t} = \{ f_s =  x^4 + sz^4 + z^2y + z = t \}$.
This deforms the $\cB_0$-type hypersurface $X_{0,t}$ into a $\cF$-type one
$X_{s,t}$ for $s\not=0$, see
 \cite[Example 6.5]{STexch}.
We recall that, for all $s$, $f_s$ has  a generic fibre,
which is homotopy equivalent to a bouquet of three 2-spheres.
There are no affine critical points and  $t=0$ is the only atypical value of $f_s$.

In order to compute the total curvature we use formulas (\ref{eq:F-form}) and (\ref{eq:B-form}) for the $\cB_0$-type ($s=0$) and
(\ref{eq:F-form}) for the $\cF$-type ($s \ne 0$). 
The input for the formulas is in the table below. The computation of
$\chi(X_{s,t})$ is via the curvature by using the Gauss-Bonnet defect.
Let us recall a few facts from \cite{STexch}:

\item{(1).} $\bar X_{s,t}$ has isolated singularities at infinity in $p:= ([0:1:0],0)$ for all $s$
 and in $q:= ([1:0:0],0)$ for $s=0$. The $\mu$'s are listed in the table.

\item{(2).} The singularities of
$\bar X_{s,t} \cap H^\ity \subset \bP^2$ change from
a single smooth line $\{ x^4=0\}$ into the
isolated point $p$ with $\tilde E_7$ singularity.

\item{(3).} The space $\bar X_{s,t}\cap \bar \cH \cap H^\ity$ has a single singularity
of type $A_3$ for $s=0$ and is smooth if $s \ne 0$.

\item{(4).} The change on the level of $\chi(\bar X_{s,t} \cap H^\ity)$ is from 2 to 5,
so $\Delta \chi^\ity =-3$. Note that $\chi^{2,4} = -4$.
In the table we use as notation $\Delta \chi = \chi^{2,4} - \chi(\bar X_{s,t} \cap H^{\ity})$.
\[
\begin{array}{|c|c|c|c|c|c|}
\hline
  & & & & & \\
(s,t) &
  \mu_p(\bar X_{s,t}) + \mu_q(\bar X_{s,t})  &
 \mu(\bar X_{s,t}\cap \bar \cH \cap H^\ity)  &
 (-1)^{n+1} \Delta \chi &
 \alpha^{(2)}_{s,t} & \chi(X_{s,t}) \\
\hline
(0,0) &  18 + 3 & 3  & 4 + 2  & 36 - 30 = 6 & 6 - 6 = 0   \\
(0,t) &  15 + 3 & 3  & 4 + 2  & 36 - 27 = 9  & 9 - 6 = 3  \\
\hline
(s,0) &  18 + 0 & -  &   9       & 36 - 27 = 9  & 9 - 9 = 0 \\
(s,t) &  15 + 0 & -  &   9       & 36 - 24 = 12 & 12 -9 = 3 \\
\hline
\end{array}
\]
NB.  In notations like $(0,t)$ we mean here that $t \ne 0$.

Let us point out that in this example we have, for each fixed $t$, a $\chi$-constant
family $X_{s,t}$ of constant degree, but with
non-constant total curvature. It turns out (by using a coordinate
change in the variable $y$) that actually this family is topologically trivial. 
 \end{example}

\subsection{Examples with non-isolated affine singularities}\label{ss:exnonisol} 

In case $Y$ is singular, one may correct the formula (\ref{eq:eulerchar})
by defining the level $n$ correction terms 
$\beta^{(n)}_Y$ as follows:
\begin{equation}\label{eq:beta3}
\chi(Y) = \chi(Y\cap \cH) + (-1)^n [\alpha^{(n)}_Y + \beta^{(n)}_Y].
\end{equation}
 Remark that we have found $\beta^{(n)}_Y$ more explicitely in case $Y$ has 
only isolated singularities, see (\ref{eq:beta2}). 
For $Y$ with non-isolated singularities,
one needs lower level corrections $\beta^{(i)}_Y$, $i\le n$, which one defines
by using equalities analogous to (\ref{eq:beta3})
for successive slices.  It follows that $\beta^{(n-i)}_Y =0$ for $i>\dim \Sing
Y$.  We show in the following
Examples \ref{e:3}  and \ref{e:4} how the lower $\beta$'s occur in case of $Y$ has
one-dimensional singularities.

\begin{example}\label{e:3} 
Consider the family given by a single polynomial $X_s= \{
f =  x^2 + x^3y + z^4 = s\}$ and note that $f$  has a non-isolated singularity.
The critical set is the $y$-axis, with constant transversal type $A_3$, and
 the only atypical value turns out to be 0.
A generic affine pencil produces a polar curve, which has $12$ intersections points
with $X_s$ if $s \ne 0$. It has $6$ intersections with $(X_{0})_\reg$ and no
intersection with $\Sing X_0$, therefore six points disappear at
infinity.
This gives the values of $\alpha^{(2)}$ in the table below.
We have $\beta^{(2)}_{X_s}$ = 0 for all $s$ since $X_s$ is non-singular
for $s\not=0$ and $X_0$ has a non-singular 1-dimensional singular locus with 
constant transversal type. 

We consider next the restriction of $f$ to a generic hyperplane section.
We use the plane  $\cH$ defined by $y = px+qz+r$.
This gives us the polynomial
\[
g = x^2 + px^4 +qx^3z + z^4 + r x^2z = s
\]
The direct computation of $\alpha^{(1)}$ turns out to be involved, so we choose
the following way. For generic $(p,q,r)$, the fibers of $g$
are general at infinity, of degree $4$. So $\chi(X_s\cap \cH) = -8$ for $s \ne
0$.
If  $s=0$ then $X_0\cap \cH$ has a $A_3$ singularity, which
has as effect $\chi(X_0\cap \cH) = -5$.
By slicing again $g$ we get $4$ points: this gives
 $\alpha^{(0)} + \beta^{(0)}$ in the table below.

Next the complex links: the fibre $X_0$ has a singular stratum which is linear
and with transversal $A_3$ singularity.
Its complex link contributes with $\beta^{(1)}=1$. If $s \ne 0$ the fibre is
smooth, so all betas are zero. Using the notations $\chi^2 = \chi(X_s)$, $\chi^1 = \chi(X_s \cap \cH)$,
$\chi^0 = \chi(X_s \cap \cH \cap \cH')$, 
the table with all information looks as follows:

\[
\begin{array}{|c|c|c|c||c||c|c|c|c|}
\hline
 \alpha^{(i)} &\beta^{(i)} &\alpha^{(i)}+\beta^{(i)} & \chi^i &i
 &\alpha^{(i)} &\beta^{(i)} &\alpha^{(i)}+\beta^{(i)} & \chi^i\\
\hline
 6 & 0 & 6 & 1 & 2 & 12 & 0 & 12 & 4 \\
 8 & 1 & 9 &-5 & 1 & 12 & 0 & 12 &-8 \\
 4 & 0 & 4 & 4 & 0 &  4 & 0 &  4 & 4 \\
\hline
\multicolumn{4}{|c||}{s = 0} & & \multicolumn{4}{|c|}{s \ne 0}\\
\hline
\end{array}
\]
\end{example}
We get $\int_{X_s} | K |\d v  =  12 \omega_2 $ if $s \ne 0$ and
$\int_{X_0} | K |\d v  =  6 \omega_2$.
The vanishing of curvature is only due to the concentration at infinity:
although we have an affine  non-isolated singularity, there is no loss of total curvature
in the affine part.
Note also:\\
$\int_{X_s \cap \cH} | K |\d v  =  12 \omega_1 $ if $s \ne 0$ and
$\int_{X_0 \cap \cH} | K |\d v  =  8 \omega_1$ and that
on this level there is an affine loss of total curvature.

\begin{example}\label{e:4}
Consider $f =  x^2y + x^3y^2 + z^5 = s$. This
can be treated in the same way.  The polynomial has a  non-isolated smooth
1-dimensional critical set ($y$-axis), but with a
non-trivial complex link on the level $i=2$ (modelled on the Whitney
umbrella) and an isolated singularity on level $i=1$.
There is also an affine contribution to the loss of total curvature.
The corresponding table is:
\[
\begin{array}{|c|c|c|c||c||c|c|c|c|}
\hline
 \alpha^{(i)} &\beta^{(i)} &\alpha^{(i)}+\beta^{(i)} & \chi^i &i
 &\alpha^{(i)} &\beta^{(i)} &\alpha^{(i)}+\beta^{(i)} & \chi^i\\
\hline
10 & 2 & 12 & 1 & 2 & 32 & 0 & 32 & 17 \\
15 & 1 & 16 &-11 & 1 & 20 & 0 & 20 &-15 \\
 5 & 0 & 5 & 5  & 0 &  5 & 0 &  5 & 5 \\
\hline
\multicolumn{4}{|c||}{s = 0} & & \multicolumn{4}{|c|}{s \ne 0}\\
\hline
\end{array}
\]

\end{example}

{\sc Acknowledgements} This research started during a one month 
visit of the first named author at the USTL, Lille.  
He also wishes to
thank for support and excellent working conditions the Institut des Hautes \'Etudes
Scientifiques at B\^ures sur Yvette. 

%


\end{document}